
\input amstex
\documentstyle {amsppt}

\document
\vsize=7 in
\hsize=5.5 in
\magnification=1200

\topmatter
\title  Determinant of complexes and higher hessians
 \endtitle
\author Fernando Cukierman
\endauthor

\leftheadtext{F. CUKIERMAN}
\rightheadtext{DETERMINANT OF COMPLEXES AND HIGHER
HESSIANS}

\address {Departamento de Matematica, Ciudad Universitaria,
(1428) Buenos Aires, ARGENTINA} \endaddress
\email {fcukier\@mate.dm.uba.ar} \endemail

\abstract {Let $X \subset \Bbb P^r$
 be a smooth algebraic curve in  projective space, over an
algebraically closed field of characteristic zero.
For each $m \in \Bbb N$,
the $m$-flexes of $X$ are defined as the points where
the osculating hypersurface of degree $m$ has higher
contact than expected,
and a hypersurface $H \subset \Bbb P^r$ is called a
$m$-Hessian if it
cuts $X$ along its $m$-flexes. When $X$ is a complete intersection,
we give an expression for a (rational) $m$-Hessian as the determinant
of a complex of graded free modules naturally associated to $X$.}
\endabstract

\endtopmatter

{\eightpoint
\flushpar
{\bf Contents.} \newline
\S 1 Introduction. \newline
\S 2  Preliminaries on differential operators. \newline
\S 3  Principal parts and closed subschemes. \newline
\S 4  Resultants, duals and Div. \newline
\S 5  Hessians.}
\newline\newline\newline

\flushpar
{\bf \S 1 Introduction.}
\newline
\newline

Our purpose in this article is to extend
a well known fact
in the theory of algebraic curves, namely, that the Hessian
of a plane curve intersects the curve in its flexes.

To recall the motivating situation,
let $X \subset \Bbb P^2$ denote a smooth
algebraic plane curve, given by an equation
$F=0$ of degree $d \ge 2$, over an algebraically closed field $k$ of
characteristic zero. The Hessian of $F$ is defined as
$ H_1(F) = \text{det } H(F)$, where
$$H(F) =
\left(\dfrac{\partial^2 F}
{\partial x_i \ \partial x_j}\right)_{0 \le i,j \le 2}$$
is the Hessian matrix of $F$, and the flexes
of $X$ (or 1-flexes, as they will be called below) are the
points where the tangent line has contact three or more.
The fact refered to before is that $(H_1(F) = 0) \cap X$
coincides with the set of flexes of $X$ (and the multiplicity
of intersection equals the weight of the flex, defined below).

At each point $x \in X$ there is a unique conic, the osculating
conic, that has maximal contact with $X$. Also, there is a finite
number of points, the 2-flexes of $X$, where the osculating conic
has contact six or more. Cayley \cite{Ca} posed himself the question of
finding a 2-Hessian, that is, a covariant curve $H_2(F)$ that cuts
$X$ along its 2-flexes. He found the expression
$$
\split
H_2(F) = (12d^2-54d+57)\ h^2 \ \text{Jac}(F,h,\Omega_H)
+ (d-2)(12d-27) \ h^2 \ \text{Jac}(F,h,\Omega_F) \\
+ 40 (d-2)^2 \ h \ \text{Jac}(F,h,\Psi)
\endsplit
$$
where
$h=H_1(F)$, \ $\Omega = (D^t H(F)^{*} D) h$, \
$\Psi = (D h)^t H(F)^{*} (D h)$,\
$D=(\partial/\partial x_0,
\partial/\partial x_1,\partial/\partial x_2)$ and
$H(F)^{*}$ is the adjoint matrix of $2 \times 2$ minors of $H(F)$.
Cayley's notation
$\Omega_H$ and $\Omega_F$ stands for certain
first partial derivatives of $\Omega$, see loc. cit.

More generally, for each $m= 1, 2, 3, \dots$ one may consider
osculating curves of degree $m$ and the corresponding $m$-flexes,
thus obtaining
a constellation of distinguished points on $X$.
The problem of finding a covariant $m$-Hessian
for $m >2$ does not seem to have been solved in the
classical literature.

Another direction of generalization is to curves
in higher dimensional projective spaces.
There is a fairly simple formula for a 1-Hessian for complete
intersection curves in $\Bbb P^3$ due to Clebsch (\cite{S}, Art. 363).
In this article we will obtain some information about $m$-Hessians
for complete intersection curves in $\Bbb P^r$.

In order to consider the more general situation, we recall at this
point some terminology regarding flexes.
Let $X$ be a smooth connected complete algebraic curve over $k$
and $L$ a line bundle over $X$.
The different orders of vanishing at $x$ of sections $s \in H^0(X,L)$
constitute a finite set with $n + 1 := \text{dim}_k  H^0(X,L)$ elements,
and when written in increasing order
$$a^L_0(x)< \dots < a^L_n(x)$$
it is called the vanishing sequence of $L$ at $x$.
The number
$w^L(x) = \sum_{i=0}^n a^L_i(x) - i $ is called the weight of $L$ at $x$,
and $x$ is said to be an $L$-flex if $w^L(x) > 0$ or,
equivalently, if there exists a section $s \in H^0(X,L)$
vanishing at $x$ to order at least $\text{dim}_k  H^0(X,L)$.

The divisor of $L$-flexes is defined as the divisor
$\sum _{x \in X } {w^L(x)}\ x$; it may be given a determinantal
expression as follows: let
$$\tau: H^0(X,L) \otimes_k \Cal O_X \to \Cal P^n_{X/k}(L)$$
denote the homomorphism obtained by
 composing the natural homomorphism \newline
$H^0(X,L) \otimes_k \Cal O_X \to L$ and the universal differential
operator $L \to \Cal P^n_{X/k}(L)$ with values in
the rank $n+1$ bundle of principal parts of order $n$ of $L$ (see \S 2).
It may be shown (see \cite{P}, \cite{ACGH})
that the divisor of $L$-flexes is equal to the divisor of zeros of
$\wedge^{n+1}\tau$.
If we choose an isomorphism
$\wedge^{n+1}H^0(X,L) \cong k$ then $\wedge^{n+1}\tau$
corresponds to a section
$w_L  \in H^0(X, \wedge^{n+1} \Cal P^n_{X/k}(L))$ that will be
called the Wronskian of $L$.

Now suppose that our curve $X$ sits in a projective space
$\Bbb P^r_k$. If $H \subset \Bbb P^r_k$ is a hypersurface, we say
that $H$ is a  $L$-Hessian if $H$ cuts $X$ along the divisor
of $L$-flexes. Notice that a necessary condition for an $L$-Hessian
to exist is that $\wedge^{n+1} \Cal P^n_{X/k}(L) \cong \Cal O_X(N)$
for some $N \in \Bbb N$. Assuming that such an isomorphism exists
and has been
chosen, an $L$-Hessian is more precisely defined as a hypersurface $H$
of degree $N$ that restricts on $X$ to $w_L$.

If $f:  X \to  S$ is a family of smooth curves in
$\Bbb P^r$ and $L$ is
a line bundle on $X$ such that $f_*(L)$ is locally free of rank $n+1$,
we define a relative $L$-Hessian to be
an effective Cartier divisor
$\Cal H \subset \Bbb P^r_{S}$
that restricts on $X$ to the relative $L$-Wronskian
$w_{f,L}  \in  H^0(X, \wedge^{n+1} \Cal P^n_{X/S}(L)
\otimes (\wedge^{n+1}f^*f_*(L))^{-1})$,
obtained from the natural homomorphism
$$\tau_f: f^*f_*(L) \to \Cal P^n_{X/S}(L)$$
Our general goal is then to write down
 an $L$-Hessian in case the
family is the universal family parametrized by a Hilbert scheme,
for a certain choice of $L$. The Hessian is also required to be
invariant under the natural action of the projective group, a  condition
that implies that one has a Hessian on each family of curves of fixed
degree and genus, compatible with pull-back.

We shall concentrate on the case of the universal family
of complete intersection curves of multi-degree
$(d_1, \dots, d_{r-1})$ in $\Bbb P^r$, choosing  $L = \Cal O(m)$
for a certain $m \in \Bbb N$.
The first remark about this case is that
an $m$-Hessian for
an individual curve
$X$ (of degree at least 3) exists: since
$\Omega^1_{X/k} \cong \Cal O_X(d)$ for
some $d \in \Bbb N$ it easily follows that the Wronskian
$w_L$ is a section of
$\Cal O_X(N)$ for some $N \in \Bbb N$, and since $X$
is projectively normal $w_L$ lifts to $\Bbb P^r_k$.

For this family, our goal amounts to the construction of
 a polynomial
$\text{GL}_{r+1}(k)$ equivariant map
$$H_m:  \prod_{j=1}^{r-1} \text{Sym}^{d_j}(k^{r+1}) \to
\text{Sym}^{N}(k^{r+1})$$
such that for each smooth complete intersection curve
$X = (F_1 = \dots = F_{r-1} = 0)$ the hypersurface
$H_m(F_1, \dots , F_{r-1})$ is an $m$-Hessian for $X$.
In (5.11) we give a  proof of existence and essential uniqueness of
such a map $H_m$, also computing its degree in each variable.
In (5.19) we give a constructive description of a rational
$m$-Hessian.

The formulas of Hesse, Cayley and Clebsch are given as the determinant
of an explicit {\it matrix} of differential operators applied to
the $F_i$.
The aim of this article is to express
$m$-Hessians
as the  Div  (see \cite{KM}) of a {\it complex}  of
graded free modules. The differentials in this complex
would be explicit matrices involving derivatives
of the polynomials $F_i$. This aim is partially achieved here:
we are able to obtain by this method only a {\it rational} $m$-Hessian
(see (5.19)).

The idea of using determinants of complexes for
writing "almost explicit" formulas is borrowed from articles by
 I. Gelfand, M. Kapranov and A. Zelevinsky; these authors took such an
approach to give
an equation for the dual hypersurface of a projective variety.

To explain our construction, let us first observe that lifting the
Wronskian $w_L$ to projective space would be immediate if the
homomorphism $\tau: H^0(X,L) \otimes_k \Cal O_X \to \Cal P^n_{X/k}(L)$
was the restriction to $X$ of a homomorphism of bundles on projective
space. Since this does not appear to be the case, we construct instead
a homomorphism $\tilde \tau$ of complexes of graded free modules
on projective space such that $\tilde \tau|_X$ is quasi-isomorphic
to $\tau$. The expression $\text{Div}(\tilde\tau)$ is then a determinantal
description of a (rational, see (5.19)) Hessian.

Since our curve is a complete intersection, the Koszul complex
provides a lifting of $H^0(X,L) \otimes_k \Cal O_X$.
On the other hand, the lifting of $\Cal P^n_{X/k}(L)$ is done
in two steps: first we use a resolution, described in \S 3,
that involves the sheaves
$\Cal P^j_{\Bbb P^r/k}(\Cal O(m))$;
these $\Cal O_{\Bbb P^r}$-Modules are not graded free
(i.e. sheafifications of graded free modules), and are then
replaced by graded free $\Cal O_{\Bbb P^r}$-Modules by means of
the Euler sequences (2.21) and a cone construction.

A related problem, considered in \cite{S},  \cite{E 1, 2}  and
\cite{Dy 1, 2},
is the one of giving equations for Gaussian maps associated
to osculating hypersurfaces of degree $m$.
We plan to show in a future article how the present set up may be used
to treat this question.
\newline\newline

It is a pleasure to thank  E. Bifet, J. Burgos, R. Piene, M. Reid and
D. Ulmer for valuable conversations and the Max-Planck-Institut f\"ur
Mathematik for its financial support and very positive mathematical
atmosphere. We are also grateful to the Fundacion Antorchas for
financial support during the final stages of this work.

\newpage

\flushpar
{\bf \S 2 Preliminaries on differential operators.}
 \newline
\newline
For ease of reference, we first review some well known
definitions and facts about differential operators.
Next we consider the graded situation,
and describe the structure
of the modules of principal parts on projective spaces (see (2.25) and
(2.31)). The algebra of differential operators on projective space is
isomorphic to the degree zero subring of the algebra
of differential operators on affine space, modulo
the two-sided ideal generated by the Euler operator, a known fact \cite{BB}
reproved in (2.23) (case $d=m=0, \ n \to \infty$).
\newline\newline
(2.1) Let $k$ be a commutative unitary ring, $A$ a commutative
$k$-algebra,
$M$ and $N$ two $A$-modules and $n$ a natural number.
A $k$-linear map $D: M \to N$ is
called an $A$-differential operator of order $\le n$
 if for all
$a \in A$ the commutator $[D,a]: M \to N$, defined by
$[D,a](m) = D(am) - a D(m)$,
is a differential operator of order $\le n-1$. Differential
operators of order $\le 0$ are, by definition, the $A$-linear
functions. The set of all such differential operators of order
$\le n$ from $M$ to $N$ will be denoted by
$$\text{D}^n_{A/k}(M,N)$$
\newline
(2.2) The $k$-module $\text{Hom}_k(M,N)$ has
two commuting structures of $A$-module, defined
by the formulas $(aL)(m)=a L(m)$ and $(La)(m)=L(am)$,
i.e. it has a structure of $(A,A)$-bimodule. The subspaces
$\text{D}^n_{A/k}(M,N)$ are stable, and hence they are
$(A,A)$-bimodules.
\newline\newline
(2.3) Let
$\text{D}^{\infty}_{A/k}(M,N)=
\bigcup_{n=0}^{\infty} \text{D}^n_{A/k}(M,N) \subset \text{Hom}_k(M,N)$
be the $(A,A)$-bimodule
of differential operators of finite order.
For three $A$-modules, composition induces a $k$-bilinear map
$$\text{D}^n_{A/k}(M_1,M_2) \times \text{D}^m_{A/k}(M_2,M_3) \to
\text{D}^{n+m}_{A/k}(M_1,M_3)$$
(this may be shown by induction, using the formula
$[D_1 \circ D_2, a] = D_1 \circ [D_2, a] + [D_1, a] \circ D_2$).
In particular,  composition defines on $\text{D}^{\infty}_{A/k}(A,A)$
a  ring structure, and hence a structure of
$(A,A)$-bialgebra.
It follows from the Jacobi formula
$[[D_1, D_2], a] = [[D_1, a], D_2] - [[D_2, a], D_1]$ and induction
that the commutator of two operators
of respective orders $\le n$ and $\le m$, has order $\le m+n-1$,
that is, the graded algebra associated to the
filtration by order is commutative.
\newline\newline
(2.4) To simplify the notation,
we will sometimes abbreviate $\text{D}^n_{A/k}(A,A) = \text{D}^n_{A/k}$.
\newline\newline
(2.5) For fixed $M$ and variable $N$, the covariant functor
$\text{D}^n_{A/k}(M,N)$ is representable:
The two $k$-algebra homomorphisms
$c, d : A \to A \otimes_k A$ defined by
$c(a) = a \otimes 1$ and $d(a) = 1 \otimes a$
give $A \otimes_k A$ a structure of $(A,A)$-bialgebra.
The kernel $I = I_{A/k}$  of the multiplication
map $A \otimes_k A  \to  A$  is a bi-ideal and
$A \otimes_k A / I^{n+1}$ has a natural structure of
$(A,A)$-bialgebra. For an $A$-module $M$  let us define
the $A$-module
$$\text{P}^n_{A/k}(M) = (A \otimes_k A / I^{n+1}) \otimes_A M$$
where, while taking $\otimes_A$ we view
$A \otimes_k A / I^{n+1}$ as an $A$-module via $d$, and
we provide $\text{P}^n_{A/k}(M)$ with an $A$-module structure
via $c$ (i.e. tensor product in the category of bi-modules). Then
$$d^n_{A/k,M}: M \to \text{P}^n_{A/k}(M)$$
defined by $m \mapsto (1\otimes_k 1)\otimes_A m$ is a differential
operator of order $\le n$, and it is universal, that is, for any
$N$ the map
$$\text{Hom}_A(\text{P}^n_{A/k}(M),N) \to \text{D}^n_{A/k}(M,N)$$
defined by composing with $d^n_{A/k,M}$, is bijective.
\newline
\newline
(2.6) In case $M=A$, instead of $\text{P}^n_{A/k}(A)$
we shall sometimes write $\text{P}^n_{A/k}$.
The universal differential operator
$d^n_{A/k}: A \to \text{P}^n_{A/k}$ takes the form
$d^n_{A/k}(a) = 1\otimes_k a$.
\newline
\newline
(2.7) For a scheme $X$ over a scheme $S$, with two
sheaves of $\Cal O_X$-modules $ M$ and $ N$,
we denote by $\Cal D^n_{X/S}(M, N)$  the sheaf of local
differential operators of order $\le n$, and by
$$d^n_{X/S,M}:  M \to \Cal P^n_{X/S}(M)$$
the universal differential operator of order $\le n$
with source $M$.
\newline
\newline
(2.8) If $X/S$ is smooth
 then we have exact sequences of $\Cal O_X$-Modules
$$0 @>>> \text{Sym}^n(\Omega^1_{X/S}) \otimes_{\Cal O_X} M @>\imath>>
\Cal P^n_{X/S}(M) @>\pi>> \Cal P^{n-1}_{X/S}(M) @>>> 0$$
\newline
(2.9) We recall from \cite{EGA IV} (16.7.3) that if $M$
is quasi-coherent (resp. coherent) then $\Cal P^n_{X/S}(M)$
is also quasi-coherent (resp. coherent). If $X/S$ is smooth
then $\Cal P^n_{X/S}$ is locally free, hence flat, and
the functor $M \mapsto \Cal P^n_{X/S}(M)$ is exact.
The natural map
$$\Cal H\text{om}_{\Cal O_X}(M, N) \to
\Cal H\text{om}_{\Cal O_X}(\Cal P^n_{X/S}(M), \Cal P^n_{X/S}(N))$$
sending $\varphi \mapsto \Cal P^n_{X/S}(\varphi)$,
is a differential operator of order $\le n$
(and, in fact, it is universal, see loc. cit. (16.7.8.2)).
\newline
\newline
(2.10) Suppose now that our $k$-algebra $A$ is $\Bbb Z$-graded:
$A = \bigoplus_{d \in \Bbb Z}  A_d$ ($\Bbb Z$ could be replaced by a
more general semigroup, but we will not follow this line here).
We provide the $k$-algebra $A \otimes_k A$ with
the total grading associated to its natural bi-grading.
The multiplication map $A \otimes_k A \to A$ is
homogeneous (of degree 0) and hence its kernel $I=I_{A/k}$
is a homogeneous ideal. It follows that for any graded
$A$-module $M$,
the $A$-module $\text{P}^n_{A/k}(M)$ has a natural grading.
\newline
\newline
(2.11) If $M$ and $N$ are graded $A$-modules and $d \in \Bbb Z$, we denote
$$\text{Hom}_k(M,N)_d = \{L \in \text{Hom}_k(M,N) / \
L(M_m) \subset N_{m+d}, \ \forall m \}$$
and
$$\text{D}^n_{A/k}(M,N)_d = \text{D}^n_{A/k}(M,N) \cap \text{Hom}_k(M,N)_d$$
If $D \in \text{D}^n_{A/k}(M,N)_d$ then we say that $D$ has order
$\le n$ and degree $d$.
\newline
\newline
(2.12) PROPOSITION. If $A$ is a finitely generated $k$-algebra and
$M$ a finitely generated $A$-module then for any $N$
$$\text{D}^n_{A/k}(M,N) = \bigoplus_{d \in \Bbb Z}  \text{D}^n_{A/k}(M,N)_d$$
that is, $\text{D}^n_{A/k}(M,N)$ is naturally a graded $A$-module.
\newline
\newline
PROOF: Let us first observe that the natural isomorphism \newline
$\alpha: \text{Hom}_A(\text{P}^n_{A/k}(M),N) \to \text{D}^n_{A/k}(M,N)$  \
takes \
$\text{Hom}_A(\text{P}^n_{A/k}(M),N)_d$ \ into \ $\text{D}^n_{A/k}(M,N)_d$ \
for any $d$.
Our hypothesis imply, by (2.8), that the $A$-module
$\text{P}^n_{A/k}(M)$ is finitely generated, and then \cite{Bo,1} (\S 21, 6,
Remark)
says that the equality
$\text{Hom}_A(\text{P}^n_{A/k}(M),N) = \bigoplus_{d \in \Bbb Z}
\text{Hom}_A(\text{P}^n_{A/k}(M),N)_d$ holds true; the proposition then
follows by applying $\alpha$.
\newline
\newline
(2.13) EXAMPLE: for a graded $A$-module $M$,
define the Euler operator
$E_M \in \text{Hom}_k(M, M)$ by
$E_M(m)=d.m$ for $m \in M_d$.
The equality $[E_M, a] = E_A(a)$ for  $a \in A$, shows that $E_M$ is a
differential operator of order $\le 1$, and degree 0.
\newline
\newline
(2.14) Now suppose $A = k[x_0, \dots, x_r]$ is a polynomial
algebra, with the $x_i$ of degree one.
We assume that the ring $k$ contains the rational numbers.
It is well known (\cite{EGA IV} (16.11) or \cite{B})
that in this case $\text{D}^{\infty}_{A/k}(A,A)$ is
isomorphic to the Weyl algebra
$k[x_0, \dots, x_r, \partial_0, \dots, \partial_r]$,
with relations $[\partial_i, x_j] = \delta_{ij}$,
$[\partial_i, \partial_j]=[x_i, x_j]=0$. Every
$D \in \text{D}^{\infty}_{A/k}(A,A)$ may be written, in a unique way,
as a finite sum
$$D = \sum_{\alpha \in \Bbb N^{r+1}} F_{\alpha}  \ \partial^{\alpha}$$
where $F_{\alpha} \in A$ and $\partial^{\alpha}$ is the differential
operator of order $\le |{\alpha}|$ and degree $- |{\alpha}|$
defined by
${\alpha}! \ \partial^{\alpha} =
\prod_{0 \le j \le r} \partial_j^{\alpha_j}$.
In the expression above, $D$ is homogeneous of degree
$d$ (see  (2.11))
iff all the $F_{\alpha}$ that appear are homogeneous
and $\text{deg}(F_{\alpha})  - |{\alpha}| =d$ for all $\alpha$.
\newline
\newline
(2.15) One may reformulate (2.14) as follows
(\cite{EGA IV} (16.4.7)):
The polynomial algebra $A[y] = A[y_0, \dots, y_r]=k[x,y]$ has a
structure of $(A,A)$-bialgebra defined by the ring homomorphisms
$c'(F(x)) = F(x)$ and $d'(F(x)) = F(x+y)$ for $F(x) \in A=k[x]$.
It is easy to check that
$$\sigma: A[y] \to A \otimes_k A$$
defined by $\sigma(x_i) = x_i \otimes 1$ and
$\sigma(y_i) = 1 \otimes x_i - x_i \otimes 1$, is an isomorphism
of $(A,A)$-bialgebras, and
since $\sigma (y_0, \dots, y_r) = I_{A/k}$, there are isomorphisms
$$\sigma: A[y]/(y)^{n+1} \to \text{P}^n_{A/k}$$
It then follows from (2.5) that the differential operator of order $\le n$
$${d'}^n_{A/k}: A \to A[y]/(y)^{n+1}$$
defined by
$${d'}^n_{A/k}(F(x)) = F(x+y) = \sum_{\alpha}
\partial^{\alpha}(F(x)) \  y^{\alpha} \ \ \ \text{mod } (y)^{n+1}$$
is universal.
\newline
\newline
(2.16) As an $A$-module via $c'$,
$A[y]/(y)^{n+1}$
is free with basis $y^{\alpha}$ for $|\alpha| \le n$. The grading
of $\text{P}^n_{A/k}$  transports  to
the grading in $A[y]/(y)^{n+1}$ such that $\text{deg}(y_i) =1$,
and the choice of basis defines an isomorphism of graded $A$-modules
$$A[y]/(y)^{n+1}  \cong \bigoplus_{|\alpha| \le n} A(-|\alpha|)
\cong \bigoplus_{j=0}^n A(-j)^{\binom {r+j}j}$$
where $M(j)$ denotes the shift by $j \in \Bbb Z$
 for a graded $A$-module $M$.
\newline
\newline
(2.17) According to (2.10), the elements
$dx_i \in \Omega^1_{A/k}$ are homogeneous of degree one. Therefore
there is an isomorphism of graded $A$-modules
$$\Omega^1_{A/k} \cong A(-1)^{r+1}$$
sending $\sum_i  a_i  dx_i \mapsto (a_i)$. The isomorphism
(2.16) is then written more intrinsically as
$$\text{P}^n_{A/k} \cong
\bigoplus_{j=0}^n  \text{Sym}^j \Omega^1_{A/k}$$
In particular, the exact sequences (2.8) split on affine space.
\newline
\newline
(2.18) (jets and polars) For $F \in A= k[x]$ define
its $j$-th polar
$$p^j(F) = \sum_{|\alpha| = j}
\partial^{\alpha}(F) \  (dx)^{\alpha} \ \
\in \text{Sym}^j \Omega^1_{A/k}$$
It follows from (2.15) and (2.17) that
the differential operator of order $\le n$
$$p^{n}_{A/k}: A \to
\bigoplus_{j=0}^n \text{Sym}^j \Omega^1_{A/k} $$
$$p^{n}_{A/k}(F) = (p^0(F), \dots, p^n(F))$$
is universal. Notice that, according to (2.11), $p^{n}_{A/k}$
is homogeneous of degree zero.
\newline
\newline
(2.19) Applying $\text{Hom}_A(-, A) = -^*$ in (2.17) we obtain
an isomorphism of graded $A$-modules
$$\text{D}^n_{A/k}(A, A) \cong
\bigoplus_{j=0}^n {\text{Sym}^j \Omega^1_{A/k}}^*$$
which corresponds to the unique representation of a
differential operator described in (2.14).
\newline
\newline
(2.20) Now we consider $P = \Bbb P^r_k = \text{Proj}(A)$.
If $M$ is a graded $A$-module, $M \sptilde$ denotes the
corresponding quasi-coherent sheaf of modules
on projective space,
$M(d)$ denotes the shifted graded module $M(d)_n = M_{n+d}$
and $\Cal O_P(m) = \Cal O(m) := A(m)\sptilde$.
\newline
\newline
(2.21) PROPOSITION (Euler sequences): For each $n \in \Bbb N$ and
$m, d \in \Bbb Z$ there is an exact sequence
of $\Cal O_P$-Modules
$$
0 @>>> \text{D}^{n-1}_{A/k}(A, A)(d)\sptilde @>.(E-m)>>
\text{D}^{n}_{A/k}(A, A)(d)\sptilde @>\rho>>
\Cal D^n_{P/k}(\Cal O(m), \Cal O(m+d)) @>>> 0
$$
where  $ E=\sum_i x_i \ \partial_i $  is the Euler differential
operator (2.13) and $.(E-m)$
is right multiplication by $(E-m) \in \text{D}^{1}_{A/k}(A, A)_0$.
\newline
\newline
PROOF: To define $\rho$ let us first remark that
for each multiplicative set $S \subset A$
there is a natural map
$$\rho_S: \text{D}^{n}_{A/k}(A, A) \to
\text{D}^{n}_{S^{-1}A/k}(S^{-1}A, S^{-1}A)$$
 which may be defined
first on derivations $\delta \in \text{Der}_{A/k}(A,A)$
by the usual quotient rule, and then using the fact that
since $A$ is a polynomial algebra, $\text{D}^{n}_{A/k}(A, A)$
is generated by the derivations. Notice
that $\rho_S$ preserves degree. Now consider
$D \in \text{D}^{n}_{A/k}(A, A)_d$,  $f \in A_{m+k}$ and $g \in A_k$,
so that $f/g$ is a section of $\Cal O_P(m)$ over the
open set where $g \ne 0$. Then $\rho(D)(f/g)$ is a section of
$\Cal O_P(m+d)$ over the same open set. Since these open sets
form a basis of open sets on $P$,  this defines the sheaf
homomorphism $\rho$.

It follows from the Euler formula that
$\rho \circ .(E-m) = 0$, and hence we have a complex that we
shall denote $(2.21)_n$. We will prove by induction on $n$ that
this complex is exact, starting with $(2.21)_0$ which is
clearly exact ($\text{D}^{-1}_{A/k} = 0$).
\newline
\newline
$(2.22)_n$ Combining (2.21) with the duals of (2.8), we obtain a
commutative diagram with exact columns,
where the maps denoted $p$ are canonical projections
$$
\minCDarrowwidth{.4 cm}
\eightpoint
\CD
  @. 0 @. 0 @. 0 @. \\
 @. @VVV  @VVV  @VVV @. \\
0 @>>> \text{D}^{n-2}_{A/k}(d)\sptilde @>.(E-m)>>
\text{D}^{n-1}_{A/k}(d)\sptilde @>>>
\Cal D^{n-1}_{P/k}(\Cal O(m), \Cal O(m+d)) @>>> 0 \\
@. @VVV  @VVV  @VVV @. \\
0 @>>> \text{D}^{n-1}_{A/k}(d)\sptilde @>.(E-m)>>
\text{D}^{n}_{A/k}(d)\sptilde @>>>
\Cal D^n_{P/k}(\Cal O(m), \Cal O(m+d)) @>>> 0 \\
@. @VpVV  @VpVV  @VVV @. \\
0 @>>>
\Cal O(d) \otimes {\text{Sym}^{n-1}(\Omega^1_{A/k})^*}\sptilde @>e>>
\Cal O(d) \otimes {\text{Sym}^{n}(\Omega^1_{A/k})^*}\sptilde @>>>
 \Cal O(d) \otimes \text{Sym}^n(\Omega^1_{P/k})^*  @>>> 0\\
@. @VVV  @VVV  @VVV @. \\
@. 0 @. 0 @. 0 @.
\endCD
$$
Here $e$ is defined as the composition
$$
\eightpoint
\Cal O(d) \otimes {\text{Sym}^{n-1}(\Omega^1_{A/k})^*}\sptilde @>.E>>
\Cal O(d) \otimes {\text{Sym}^{n-1}(\Omega^1_{A/k})^*}\sptilde
\  \oplus \
\Cal O(d) \otimes {\text{Sym}^{n}(\Omega^1_{A/k})^*}\sptilde @>p>>
\Cal O(d) \otimes {\text{Sym}^{n}(\Omega^1_{A/k})^*}\sptilde
$$
and is given by the formula (see (2.24))
$$e(\sum_{|\alpha|=n-1} F_{\alpha} \ \partial^{\alpha}) =
\sum_{\alpha, i} \ x_i \ F_{\alpha} \ \partial_i \ \partial^{\alpha}$$
In other words, $e$ is multiplication by $E$ in the symmetric
algebra of the module $({\Omega^1_{A/k}}^*)\sptilde$. The lower row
of $(2.22)_1$ is dual to the usual Euler sequence
$$0 @>>> \Omega^1_{P/k} @>{\jmath}>> {\Omega^1_{A/k}}\sptilde
@>\varepsilon>> \Cal O_P @>>> 0$$
which is known to be exact ($\varepsilon$ is contraction with $E$ and
$\jmath$ is pull-back of 1-forms to affine space). Then for any
$n$ the lower row of $(2.22)_n$ is  exact, since it is obtained from
the usual Euler sequence by applying $\text{Sym}^n$
(see \cite{I}, (4.3.1.7)). Now suppose
$(2.21)_{n-1}$ is exact. Considering the long exact sequence of
homology of $(2.22)_n$ it follows that $(2.21)_n$ is also exact, which
finishes the proof.
\newline
\newline
(2.23) COROLLARY: there are exact sequences
$$
0 @>>> \text{D}^{n-1}_{A/k}(A, A)_d @>H^0(.(E-m))>>
\text{D}^{n}_{A/k}(A, A)_d @>>>
H^0(P,\Cal D^n_{P/k}(\Cal O(m), \Cal O(m+d))) @>>> 0
$$
\newline
PROOF: We take cohomology in (2.21) and observe that
 (2.19) and the calculation of cohomology of line bundles on $P$
imply that
$H^0(P, \text{D}^{n}_{A/k}(A, A)(d)\sptilde) = \text{D}^{n}_{A/k}(A, A)_d$
and $H^1(P,\text{D}^{n}_{A/k}(A, A)(d)\sptilde) = 0$
for all $n \in \Bbb N$ and  $d \in \Bbb Z$.
\newline
\newline
(2.24)  To make the map $H^0(.(E-m))$ explicit, in view of (2.19),
we need to write down the standard form (2.14) for $D.(E-m)$
for each given
$D = \sum_{\alpha} F_{\alpha}\ \partial^{\alpha} \in
\text{D}^{n-1}_{A/k}(A,A)_d$. This is achieved by combining
$D.(E-m) = D.E - mD$ with the easily checked formula
$$D.E = \sum_{\alpha}\ (\sum_i \ x_i \ F_{\alpha}  \ \partial_i
\ \partial^{\alpha})  +  |{\alpha}| \ F_{\alpha}  \ \partial^{\alpha}$$
\newline
(2.25) Applying $\Cal H\text{om}_{\Cal O_P}( - , \Cal O(m+d))$
in (2.21) we obtain the following exact
sequence that describes
the structure of $\Cal P^n_{P/k}(\Cal O(m))$
$$
0 @>>> \Cal P^n_{P/k} (\Cal O(m)) @>\jmath>>
\text{P}^n_{A/k}(A)(m)\sptilde @>.(E-m)^*>>
\text{P}^{n-1}_{A/k}(A)(m)\sptilde @>>> 0
$$
\newline
(2.26) Combining (2.25) with (2.17) and (2.18) we obtain a
commutative diagram that describes the universal
differential operator of order $\le n$ for $\Cal O(m)$
$$
\CD
0 \to \Cal P^n_{P/k} (\Cal O(m)) @>>>
\displaystyle
\bigoplus_{j=0}^n  (\text{Sym}^j \Omega^1_{A/k})(m)\sptilde
@>.(E-m)^*>>
\displaystyle
\bigoplus_{j=0}^{n-1}  (\text{Sym}^j \Omega^1_{A/k})(m)\sptilde  \to 0 \\
  @A{d^n_{P/k,\Cal O(m)}}AA  @Ap^{n}_{A/k}AA @. @. \\
\Cal O(m)  @= \Cal O(m)
\endCD
$$
\newline
(2.27) Let us make $.(E-m)^*$ explicit. It follows from (2.24)
that $.(E-m)^*$ is given on basis elements by
$$.(E-m)^*((dx)^\beta) =  c.(dx)^\beta +
\sum_ {\beta_i \ge 1} x_i (dx)^{\beta-e_i}$$
where $e_i \in \Bbb N^{r+1}$, \ $(e_i)_j = \delta_{ij}$
and $c= |\beta| - m$ for $|\beta| < n$,\ $c=0$  for $|\beta| = n$.
\newline
In other words, $.(E-m)^*$ is the direct sum of maps
$$(c_j, \varepsilon): (\text{Sym}^j \Omega^1_{A/k})(m)\sptilde \to
(\text{Sym}^j \Omega^1_{A/k})(m)\sptilde \ \oplus \
(\text{Sym}^{j-1} \Omega^1_{A/k})(m)\sptilde$$
where $c_j$ is multiplication by $j - m$ for $j < n$,\ $c_n=0$,\  and
$\varepsilon$ is contraction with $E$.
\newline
\newline
Now we take a closer look at (2.25)-(2.27); we shall obtain
a better description of the modules $\Cal P^n_{P/k} (\Cal O(m))$.
\newline
\newline
$(2.28)_n$ We start by writing down the dual of $(2.22)_n$
$$
\minCDarrowwidth{.4 cm}
\eightpoint
\CD
  @. 0 @. 0 @. 0 @. \\
 @. @AAA  @AAA  @AAA @. \\
0 @>>> \Cal P^{n-1}_{P/k}(\Cal O(m)) @>\jmath>>
\text{P}^{n-1}_{A/k}(m)\sptilde @>.(E-m)^*>>
\text{P}^{n-2}_{A/k}(m)\sptilde @>>> 0 \\
@. @A{\pi}AA  @A{\pi}AA  @A{\pi}AA @. \\
0 @>>>
\Cal P^n_{P/k}(\Cal O(m))  @>\jmath>>
\text{P}^{n}_{A/k}(m)\sptilde @>.(E-m)^*>>
\text{P}^{n-1}_{A/k}(m)\sptilde @>>> 0 \\
@. @A{\imath}AA  @A{\imath}AA  @A{\imath}AA @. \\
 0 @>>>
\text{Sym}^n \Omega^1_{P/k} \otimes \Cal O(m)  @>{\jmath}>>
(\text{Sym}^{n} \Omega^1_{A/k})(m)\sptilde
@>\varepsilon>>
(\text{Sym}^{n-1} \Omega^1_{A/k})(m)\sptilde @>>> 0\\
@. @AAA  @AAA  @AAA @. \\
@. 0 @. 0 @. 0 @.
\endCD
$$
\newline
(2.29) The composition of $\jmath$ with the projection
$p:  {\text{P}^n_{A/k}(m)}\sptilde \cong
\bigoplus_{j=0}^n  (\text{Sym}^j \Omega^1_{A/k})(m)\sptilde \to
(\text{Sym}^n \Omega^1_{A/k})(m)\sptilde$ will be denoted
$$\varphi= \varphi_{m,n}: \Cal P^n_{P/k}(\Cal O(m)) \to
(\text{Sym}^n \Omega^1_{A/k})(m)\sptilde$$
It is a homomorphism between bundles of the same rank. We are
interested in the bundle on the left, and since the one on the right
is understood (it is graded free), we shall analyse to what extent
$\varphi$
is an isomorphism. The answer depends on $m$ and $n$ and is given in
the next Theorem.
\newline
\newline
(2.30) We need another piece of notation, namely
$$\lambda_j(m,n) := \prod_{i=n-j}^{n-1} (m-i)^{-1} \ \ \in \Bbb
Q$$ which is well defined as long as $i \ne m$ for
$n-j \le i \le n-1$, that is, $m < n-j$ or $m \ge n$.
When the values of $m$ and $n$ are understood, we will write
$\lambda_j$ instead of $\lambda_j(m,n)$.
\newpage

\flushpar
(2.31) THEOREM. With the notation introduced above, the following
is true:
\newline
\newline
(a) If $m < 0$ or if $m \ge n$ then $\varphi$ is an isomorphism.
Its inverse is given by the formula
$$\varphi^{-1}(w_n) = \sum_{j=0}^n \ \ \lambda_j \
\varepsilon^j(w_n)$$ where
$$\varepsilon^j: \text{Sym}^{i} \Omega^1_{A/k} \to
\text{Sym}^{i-j} \Omega^1_{A/k}$$
is iterated contraction with the Euler operator $E$.
\newline\newline
(b) If $0 \le m < n$ then there is an exact sequence
$$
0 @>>>
\text{Sym}^m \Omega^1_{A/k}(m)\sptilde
@>t>>
\Cal P^n_{P/k}(\Cal O(m))
@>{\varphi}>>
\text{Sym}^n \Omega^1_{A/k}(m)\sptilde
@>{\varepsilon^{n-m}}>>
\text{Sym}^m \Omega^1_{A/k}(m)\sptilde
@>>> 0
$$
where
$$t(w_m) = \sum_{i=0}^m \frac{\varepsilon^i(w_m)}{i !}$$
\newline
PROOF: A typical section $w$ of
$\bigoplus_{j=0}^n  (\text{Sym}^j \Omega^1_{A/k})(m)\sptilde$
over an (unspecified) open set may be written as $w = \sum_{j=0}^n w_j$
where $w_j$ is a rational symmetric form  of degree $j$ on affine
space; more precisely, $w_j = \sum_{|\alpha|=j}\  a_{\alpha}(x)\
(dx)^{\alpha}$
where the $a_{\alpha}$ are homogeneous rational functions of degree
$m$.
According to (2.25), such a $w$ is a section of
$\Cal P^n_{P/k}(\Cal O(m))$ if and only if $.(E-m)^*(w) = 0$. Using
(2.27) we find
$$.(E-m)^*(w) = \sum_{j=0}^n \ .(E-m)^*(w_j) = \varepsilon(w_n) +
\sum_{j=0}^{n-1} \varepsilon(w_j) + (j-m) w_j =
\sum_{j=1}^{n} \varepsilon(w_j) + \sum_{j=0}^{n-1} (j-m) w_j$$
which equals zero if and only if
$$ (m-j)\  w_j = \varepsilon(w_{j+1}) \ \ \
\text{for all } j= 0, 1, \dots, n-1 \tag*$$
Notice that with the present notation
$\varphi$ is defined by $\varphi(\sum_{j=0}^n w_j) = w_n$.
\newline
\newline
Now we consider our two cases.
\newline\newline
(a) Here $m-j \ne 0$ for $0 \le j \le n-1$ and then
$w_j = (m-j)^{-1}\  \varepsilon(w_{j+1})$. It follows that
all the $w_j$ are determined by $w_n$ and it is easy to verify
the proposed formula for $\varphi^{-1}$.
\newline
\newline
(b) Since the index $j=m$ now occurs, $w_n$ does not determine
all the other $w_j$ and $\varphi$ will have a non-zero kernel and
cokernel. To determine these, we separate $(^*)$ into two groups
of equations:
$$
\aligned
m\  w_0 &= \varepsilon(w_1) \\
(m-1)\  w_1 &= \varepsilon(w_{2}) \\
\dots \\
 w_{m-1} &= \varepsilon(w_m)
\endaligned
\qquad
\text{and}
\qquad
\aligned
0 = (m-m)\  w_m &= \varepsilon(w_{m+1}) \\
(-1)\  w_{m+1} &= \varepsilon(w_{m+2}) \\
\dots \\
(m-n+1)\  w_{n-1} &= \varepsilon(w_n)
\endaligned
$$
Notice that there are no conditions on $w_m$ and that $w_m$
determines $w_j$ with $0 \le j \le m$.
\newline
Suppose that $\varphi(w) = w_n = 0$. It follows from the second
group of equations that \newline $w_{m+1}, \dots, w_n$ are all zero,
which implies that $\text{image}(t) = \text{kernel}(\varphi)$.
\newline
On the other hand, the condition $\varepsilon(w_{m+1}) = 0$
is equivalent to $\varepsilon^{n-m}(w_{n}) = 0$, and it follows
that $\text{image}(\varphi) = \text{kernel}(\varepsilon^{n-m})$.
\newline
Finally, the surjectivity of $\varepsilon^{n-m}$ follows from the
surjectivity of $\varepsilon$, proved in (2.22).
\newline
\newline
(2.32) If $X$ is a scheme over a scheme $S$ and
$\varphi: M \to N$ is a homomorphism of $\Cal O_X$-Modules,
one has an induced homomorphism
$\Cal P^n_{X/S}(\varphi): \Cal P^n_{X/S}(M) \to \Cal P^n_{X/S}(N)$
(as in (2.9)) that involves "derivatives of $\varphi$".
With a view towards our applications, now we make
this explicit, at least in the cases in which $X$ is affine space and
$M=N=\Cal O_X$, and $X$ is projective space and $M$ and $N$
are (direct sums of) line bundles.
\newline
\newline
(2.33) For a ring $R$ and element $r \in R$ let us denote
$.r: R \to R$ the operator of multiplication by $r$. Let $A$
be a $k$-algebra as in (2.5) and $a \in A$. The operator $.a$
induces a map (2.9)
$$\text{P}^n_{A/k}(.a): \text{P}^n_{A/k}(A) \to \text{P}^n_{A/k}(A)$$
which is simply the one induced by
$.(1\otimes a): A \otimes_k A \to A \otimes_k A$.
In case $A$ is a polynomial algebra, under the isomorphism
$\sigma$ of (2.15) $\text{P}^n_{A/k}(.a)$ corresponds to multiplication
by
\newline
$a(x+y)$. Using the notation of polars as in (2.18),
since $a(x+y) = \sum_{j=0}^n p^j(a)$, multiplication
by $a(x+y)$ is a direct sum of multiplication operators
in the symmetric algebra of $\Omega^1_{A/k}$
$$.p^{j-i}(a): \text{Sym}^i \Omega^1_{A/k} \to
\text{Sym}^j \Omega^1_{A/k}$$
for $0 \le i \le j \le n$.
\newline
\newline
(2.34) We remark at this point that we also have a homomorphism
$$.a': \text{P}^{n-1}_{A/k}(A) \to \text{P}^n_{A/k}(A)$$
(as in Example (3.7)) induced by multiplication by
$a' = 1\otimes a - a \otimes 1$. It follows easily from the
definitions that its dual is the homomorphism
$$[\ ,a]: \text{D}^n_{A/k}(A,A)  \to  \text{D}^{n-1}_{A/k}(A,A)$$
sending $D \mapsto [D,a] = D.a-a.D$. In case $A$ is a
polynomial algebra $.a'$ corresponds, under the isomorphism
$\sigma$ of (2.15), to
multiplication by $\Delta a = a(x+y) - a(x)$.
\newline
\newline
(2.35) Consider $A = k[x_0, \dots, x_r]$ and let $F \in A_d$ be
a homogeneous polynomial of degree $d$. Then $F$ induces
a homomorphism $.F: \Cal O_P(m) \to \Cal O_P(m + d)$ for
each $m \in \Bbb Z$,
and we wish to describe $\Cal P^n_{P/k}(.F)$.
Keeping track of gradings in (2.33), one has a homomorphism
of degree zero
$$\text{P}^n_{A/k}(.F): \text{P}^n_{A/k}(A)(m) \to \text{P}^n_{A/k}(A)(m+d)$$
Taking (2.25) into account, $\Cal P^n_{P/k}(.F)$ is described
as the restriction of $\text{P}^n_{A/k}(.F)\sptilde$ to the submodule
$\Cal P^n_{P/k} (\Cal O(m))$. Summarizing, $\Cal P^n_{P/k}(.F)$
is explicitly given as a direct sum of maps which are multiplication
by polars of $F$ as in (2.33).
\newpage

\flushpar
{\bf \S 3 Principal parts and closed subschemes.}
\newline
\newline
In this section we study some resolutions (see (3.5) and (3.9))
of the modules of principal parts of a closed subscheme,
constructed out of a given resolution of the structure sheaf
of the subscheme.
\newline
\newline
We start by recalling \cite{EGA IV} (16.4.20):
\newline
\newline
(3.1) PROPOSITION. Let $Y$ be a scheme over a scheme $S$ and
$f: X \to Y$ a closed immersion corresponding to a quasi-coherent
Ideal $ J \subset \Cal O_Y$. Then for any $n \in \Bbb N$ the
natural map
$f^* \Cal P^n_{Y/S} \to \Cal P^n_{X/S}$ is surjective, and its kernel
is the Ideal of $f^* \Cal P^n_{Y/S}$ generated by
$f^* (\Cal O_Y.d^n_{Y/S}(J))$.
\newline
\newline
(3.2) We remark that in the affine case where $S=\text{Spec}(k)$,
$Y=\text{Spec}(A)$
and $X=\text{Spec}(B)$ with $B=A/J$, (3.1) is equivalent to the fact
that the natural map
$$A \otimes_k A / (I^{n+1}_{A/k} + J \otimes_k A+A \otimes_k J)
 \to B \otimes_k B / I^{n+1}_{B/k}$$
is an isomorphism.
\newline
\newline
(3.3) It follows from (3.1) that we have an exact sequence of
$\Cal O_X$-Modules
$$
f^* \Cal P^n_{Y/S}(J)  @>{f^* \Cal P^n_{Y/S}(i)} >> f^* \Cal P^n_{Y/S}
@>>> \Cal P^n_{X/S} @>>> 0
$$
where $i$ denotes the inclusion $J \to \Cal O_Y$.
Combining this with the exactness of
the functor $\Cal P^n_{Y/S}$ (now assuming $Y/S$ smooth)
and the right-exactness of $f^*$,
we deduce that for each exact
sequence of $\Cal O_Y$-Modules
$$
M @>\mu>> \Cal O_Y @>>> \Cal O_X @>>> 0
$$
the sequence
$$
f^* \Cal P^n_{Y/S}(M)  @>f^* \Cal P^n_{Y/S}(\mu) >>
f^* \Cal P^n_{Y/S}
@>>> \Cal P^n_{X/S} @>>> 0
$$
is exact.
\newline
\newline
(3.4) Suppose that we are given a resolution $(M,\mu)$ of
$\Cal O_X$,
that is, an exact complex of quasi-coherent
$\Cal O_Y$-Modules of the form
$$
 M_r @>\mu_r>> \dots @>>> M_1 @>\mu_1>>
M_0 =\Cal O_Y @>>> \Cal O_X @>>> 0
$$
We are interested in functorially associating to $(M, \mu)$ a
resolution of $\Cal P^n_{X/S}$,
continuing the first step (3.3). More precisely, we wish
(in principle, see (3.10)) to construct a complex
$(M^n, \mu^n)$ of $\Cal O_Y$-Modules such that
$f^*(M^n, \mu^n)$ is a  resolution of $\Cal P^n_{X/S}$
by $\Cal O_X$-Modules.
\newline
\newline
(3.5) Let us consider the case in which
$S=\text{Spec}(k)$, $Y=\text{Spec}(A)$
and $X=\text{Spec}(B)$ with $B=A/J$, where the ideal $J$
is generated by $a_1, \dots, a_s \in A$.
We denote $R=A \otimes_k A$, $I = I_{A/k}$ as in (2.5) and
$a'_i = 1 \otimes a_i - a_i \otimes 1$, and consider the
Koszul complex $K(a')$ of the elements $a'_i \in R$. Writing
$V = k^s$ this complex looks like
$$
0 @>>> \wedge^s V \otimes_k R @>\alpha_s>> \dots @>>>
\wedge^2 V \otimes_k R @>\alpha_2>>  V \otimes_k R
@>\alpha_1>> R
$$
We observe that since $a'_i \in I$, the differentials $\alpha_j$
correspond to matrices with all the entries in $I$, that is,
$\alpha_j(I^m.(\wedge^j V \otimes_k R)) \subset
I^{m+1}.(\wedge^{j-1} V \otimes_k R)$ for any $m$. It follows that
there is an induced complex (interpreting $I^j=R$ for $j\le 0$)
$$
0 @>>> \wedge^s V \otimes_k R/I^{n-s+1} @>>> \dots
@>>>  \wedge^2 V \otimes_k R/I^{n-1}
@>>> V \otimes_k R/I^{n} \to R/I^{n+1}
$$
that is,  a complex
$$
0 @>>> \wedge^s V \otimes_k \text{P}^{n-s}_{A/k} @>\alpha_s>> \dots
@>>>     \wedge^2 V \otimes_k \text{P}^{n-2}_{A/k}
@>\alpha_2>>  V \otimes_k \text{P}^{n-1}_{A/k} @>\alpha_1>> \text{P}^n_{A/k}
$$
that we shall denote $\text{P}^n_{A/k}(a_1, \dots,a_s) = \text{P}^n_{A/k}(a)$.
\newline
\newline
(3.6) PROPOSITION. Suppose that $a_1, \dots, a_s \in A$ is a
regular sequence and that $A$ and $B$ are smooth over $k$.
Then $B \otimes_A \text{P}^n_{A/k}(a)$  is a resolution of
$\text{P}^n_{B/k}$.
\newline
\newline
PROOF: Denote by $\Cal P^n$ the augmented complex
$B \otimes_A \text{P}^n_{A/k}(a) \to \text{P}^n_{B/k} \to 0$. We will
show that $\Cal P^n$ is exact by induction on $n$.
The statement is obvious for $n=0$. Assuming it
for $n-1$, consider the natural surjective map
$\Cal P^n \to \Cal P^{n-1}$ and let $\Cal Q^n$ denote its kernel.
The complex $\Cal Q^n$ looks like
$$
\eightpoint
0 \to \wedge^s V \otimes_k I^{n-s}/I^{n-s+1} \to \dots
\to  V \otimes_k I^{n-1}/I^{n} \to I^{n}/I^{n+1}
\to I_B^{n}/I_B^{n+1} \to 0
$$
and it has a natural map, call it $\sigma$, from the complex
$\Cal {Q'}^n$:
$$
\eightpoint
0 \to \wedge^s V \otimes_k \text{Sym}^{n-s}(I/I^2) \to \dots
\to  V \otimes_k \text{Sym}^{n-1}(I/I^2) \to \text{Sym}^{n}(I/I^2)
\to \text{Sym}^{n}(I_B/I_B^2) \to 0
$$
Since $A/k$ and $B/k$ are smooth, all the $A$ and $B$-modules in the
complexes above are locally free.
Therefore, the exact
sequences $0 \to \Cal Q^n \to \Cal P^n \to \Cal P^{n-1} \to 0$
remain exact after $B \otimes_A -$. Using the long
exact sequence of homology, we are reduced to showing that
$B \otimes_A\Cal Q^n$ is exact for all $n \ge 1$. By the smoothness
hypothesis again, $\sigma$ is an isomorphism, and hence we only
need to show that $B \otimes_A\Cal {Q'}^n$ is exact. But this complex is
exact because it is isomorphic to the $n$-th graded piece of the
homomorphism
of Koszul algebras
$$\wedge(J/J^2) \otimes_B
\text{Sym}(B\otimes_A \Omega^1_{A/k}) \to
\text{Sym}(\Omega^1_{B/k})$$
associated to the exact sequence
$$ 0 \to J/J^2 \to
B\otimes_A \Omega^1_{A/k} \to \Omega^1_{B/k} \to 0$$
as in \cite{I} (4.3.1.7). This finishes the proof.
\newline
\newline
(3.7) EXAMPLE: for $s=1$, $B = A/A.a$ and the resolution
in (3.5) is
$$
0 @>>> B \otimes_A \text{P}^{n-1}_{A/k} @>.a'>> B \otimes_A \text{P}^n_{A/k}
 @>>> \text{P}^n_{B/k} @>>> 0
$$
where $.a'$ is multiplication by $a' =1\otimes a - a \otimes 1$ (see
(2.34)).
\newline
\newline
(3.8) EXAMPLE: In case $A=k[x]$ is a polynomial algebra,
denoting $\Delta a = a(x+y) - a(x) \ \in A[y]$ for $a \in A$,
the complex $\text{P}^n_{A/k}(a_1, \dots,a_s)$ of (3.5) is induced
by the Koszul complex of
$\Delta a_1, \dots, \Delta a_s$, and is written
$$
0 @>>> \wedge^s V \otimes A[y]/(y)^{n-s+1} @>\Delta a>> \dots
@>>>     \wedge^2 V \otimes A[y]/(y)^{n-1}
@>\Delta a>>  V \otimes A[y]/(y)^{n} @>\Delta a>>
A[y]/(y)^{n+1}
$$
\newline
Returning to the general problem in (3.4), and keeping the
notation from there, we have the following proposition
\newline
\newline
(3.9) PROPOSITION. Assume $Y/S$ is smooth and let
$(M,\mu)$  be a complex of
$\Cal O_Y$-Modules such that $\mu_j(M_j) \subset J. M_{j-1}$
(i.e. $f^*(\mu_j) = 0$) for all $j$.
Then for each $n \in \Bbb N$ there is a complex
$\Cal P^n(M,\mu)$ of $\Cal O_X$-Modules, functorial in $(M,\mu)$,
of the form
$$
 \dots    @>>> f^* \Cal P^{n-2}_{Y/S}(M_2) @>>>
f^* \Cal P^{n-1}_{Y/S}(M_1) @>>> f^* \Cal P^n_{Y/S}(M_0)
$$
where we use the convention $\Cal P^m_{Y/S} = 0$ for $m < 0$.
\newline
\newline
PROOF: Let $\varphi: M \to N$ be a homomorphism of
$\Cal O_Y$-Modules such that $f^*(\varphi) = 0$ and consider
the diagram below, obtained by applying (2.8) to $M$ and $N$
$$
\CD
0 @>>> \text{Sym}^n(\Omega^1_{Y/S}) \otimes M @>>>
 \Cal P^n_{Y/S}(M) @>>>  \Cal P^{n-1}_{Y/S}(M) @>>> 0 \\
@.  @VV{1\otimes\varphi}V  @VV{\Cal P^n_{Y/S}(\varphi)}V
@VV{\Cal P^{n-1}_{Y/S}(\varphi)}V  \\
0 @>>> \text{Sym}^n(\Omega^1_{Y/S}) \otimes N @>>>
 \Cal P^n_{Y/S}(N) @>>>  \Cal P^{n-1}_{Y/S}(N) @>>> 0
\endCD
$$
Now we apply $f^*$ and observe that
$f^*(1 \otimes \varphi)= 1 \otimes f^*(\varphi) = 0$.
Therefore, $f^*\Cal P^n_{Y/S}(\varphi)$ factors through
$f^*\Cal P^{n-1}_{Y/S}(M)$; denote
$$\Cal Q^n(\varphi): f^*\Cal P^{n-1}_{Y/S}(M) \to
f^*\Cal P^n_{Y/S}(N)$$
the induced homomorphism.

If $\psi: L \to M$ is another
homomorphism of $\Cal O_Y$-Modules such that $f^*(\psi)=0$
and such that $\varphi \circ \psi =0$
then an easy diagram chase shows that
$\Cal Q^{n+1}(\varphi) \circ \Cal Q^n(\psi) = 0$ for all $n \in \Bbb N$.

Now, given the complex $(M,\mu)$ and $n \in \Bbb N$,
the desired complex $\Cal P^n(M,\mu)$
is defined as
$$
 \dots    @>\Cal Q^{n-2}(\mu_3)>> f^* \Cal P^{n-2}_{Y/S}(M_2)
@>\Cal Q^{n-1}(\mu_2)>> f^* \Cal P^{n-1}_{Y/S}(M_1)
@>\Cal Q^{n}(\mu_1)>> f^* \Cal P^n_{Y/S}(M_0)
$$
\newline
(3.10) REMARK. An easy diagram chase shows that the natural
homomorphisms \newline
$\tau_j: H^0(Y,M_j) \otimes \Cal O_X \to f^*\Cal P^{n-j}_{Y/S}(M_j)$
commute with the differentials and hence
give a homomorphism of complexes
$$\bar\tau: H^0(Y,(M,\mu)) \otimes \Cal O_X \to \Cal P^n(M,\mu)$$
\newline
(3.11) REMARK. The complex $\Cal P^n(M,\mu)$
is not necessarily $f^*$ of a complex of $\Cal O_Y$-Modules
since the homomorphisms $\Cal Q^{i}(\mu_j)$ are in principle
defined only over $X$. However, in the special situation
of (3.5) the two complexes considered coincide and are
defined over $Y$.

\newpage

\flushpar
{\bf \S 4 Resultants, duals and Div.}
\newline
\newline
In this section we introduce some of the
ideas of \cite{GKZ} in order to motivate our construction
in \S 5.
\newline
\newline
(4.1) Let us consider a smooth irreducible projective variety
$X \subset \Bbb P(V)$ over an algebraically
closed field $k$ of characteristic zero.
The dual variety
$X^* \subset \Bbb P(V^*)$ is defined as the set of hyperplanes
$H \subset \Bbb P(V)$ such that $X \cap H$ is singular.
\newline
We denote by $j^1$  the composition
$V^* \to H^0(X, \Cal O (1)) \to
H^0(X, \Cal P^1_{X/k}(\Cal O (1)))$
where the second map is $H^0(d^1_{X/k,\Cal O (1)})$ as in (2.7).
Then $X^*$ may also be described as the set of linear forms
$s \in V^*$ such that their first jet $j^1(s)$ has a zero $x \in X$.
\newline
If $X^*$ is a hypersurface in $\Bbb P(V^*)$ (which is the case
for most $X$'s) then we denote by $\delta_X$ a defining equation
for $X^*$, so that $\delta_X(H) = 0$ iff $X \cap H$ is singular.
For example \cite{GKZ}, if $X$ is a Veronese variety
then $\delta_X$ is a discriminant, and if $X$ is a convenient
Segre variety then $\delta_X$ is a hyperdeterminant.
\newline
\newline
(4.2) Let $E$ be a vector bundle of rank $r$ on
a variety $X$ of dimension $n$, and  denote by
$R_E \subset \Bbb P H^0(X,E)$ the subvariety of sections
of $E$ that have a zero $x \in X$. It is not hard to check that
if $r=n+1$ and $E$ is generated by its global sections then the
resultant variety $R_E$ is an irreducible
hypersurface with degree equal to $\text{deg }c_n(E)$;
let us denote by $\rho_E$
an equation for $R_E$.
As examples of resultants
we may mention:
\newline
(a) If $X = \Bbb P^n$ and $E = \bigoplus_{j=0}^n \Cal O(d_i)$
then $R_E$ is the classical resultant hypersurface, that is, the
 collection of $n+1$-tuples of polynomials
$F_0, \dots, F_n$ of degrees
$d_0, \dots, d_n$ that have a common zero
$x=(x_0, \dots, x_n) \in \Bbb P^n$.
\newline
(b) If $X^n \subset \Bbb P^r$ and $E = \Cal O_X(1)^{n+1}$
then $R_E$ is closely related to the Chow form of $X$.
\newline
(c) If $X^n \subset \Bbb P^r$ and $E = \Cal P^1_{X/k}(\Cal O(1))$
then, according to the discussion in (4.1), we have
$$R_E \cap \Bbb P H^0(X, \Cal O(1)) = X^*$$
via the natural inclusion
$H^0(X, \Cal O(1)) \hookrightarrow
H^0(X, \Cal P^1_{X/k}(\Cal O(1)))$. In particular,  dual varieties
are linear sections of resultants,  and we also see that if
$\Omega^1_X (1)$,  and hence $E$, is generated
by global sections then $R_E$ and  $X^*$ are hypersurfaces.
\newline
\newline
(4.3) It seems to be of some interest to write down
explicit formulas for resultant polynomials  $\rho_E$.
In general this is not an easy
task, already for instance in the classical case (4.2)(a) for $n > 1$.
But  I. Gelfand, M. Kapranov and A. Zelevinsky have shown
that it is easy to give an expression for all resultants in terms
of a single operation, the determinant functor.
\newline
\newline
(4.4)  In the literature one finds at least two notions of determinant
for complexes: the Whitehead torsion
\cite{Bu}, \cite{Mi}, \cite{R} and the
Grothendieck determinant functor \cite{Bu}, \cite{KM}.

In one of its variations,
the Whitehead torsion is a function
$\tau_R$ that assigns to each homotopy equivalence
$\varphi: E \to F$
between two finite complexes $E$ and $F$
of finitely generated projective modules over a ring $R$, an
element $\tau_R(\varphi)$ in the group $K_1(R)$.
The function $\tau_R$ satisfies
some properties that make it, in the terminology of
\cite{Bu} (1.2), a theory of determinants with values in $K_1(R)$.
The torsion of a homotopically trivial complex is  then
defined through
the choice of a trivialization, see \cite{R}.

On the other hand,
the Grothendieck determinant functor
$\text{det}:  \Cal C^{\cdot}is_X \to \Cal Pis_X$ is a functor
from the category of bounded complexes of finite locally free
modules over a scheme $X$, with quasi-isomorphisms of
complexes as morphisms, to the category of graded invertible
modules on $X$, with module isomorphisms as morphisms.
The functor det satisfies some properties that make it,
in the terminology of \cite{Bu} (1.6), a general theory of
determinants with values in the Picard category $\Cal Pis_X$.
The existence of det is proved in \cite{Bu} for $X$ affine
and in \cite{KM} for general $X$.

The relationship between theories of determinants with
values in abelian groups and general theories of determinants
with values in Picard categories is summarized in a question
by Grothendieck \cite{Bu}, page 37.

In \cite{KM} it is also proved that det extends
to the category of perfect complexes, which means in
particular that each
perfect complex has associated to it a determinant line bundle,
a fact that is sometimes used to define certain line bundles
on moduli spaces (see \cite{BGS} for example).
The second ingredient of the functor det, namely its action on
morphisms, has special interest for us since it leads
to the Div construction: If
$\varphi: E \to F$ is a morphism in $\Cal C^{\cdot}is_X$ then
there is an induced isomorphism of line bundles
$\text{det}(\varphi): \text{det}(E) \to \text{det}(F)$, whereas
if $\varphi: E \to F$ is a (good) homomorphism of complexes,
not necessarily a quasi-isomorphism as above, then
$\text{Div}(\varphi)$ is a Cartier divisor supported on
the closed subscheme where $\varphi$ is not a
quasi-isomorphism (see \cite{KM} for  definition
and basic properties of Div).
\newline
\newline
(4.5) In more concrete terms, the determinant of an exact complex
of finite dimensional vector spaces with chosen basis
is obtained as an alternating
product of determinants of certain well chosen minors of the matrices of
the complex \cite{GKZ,2}. A complex of graded free modules
on projective space is determined by a collection of matrices of
homogeneous polynomials; if it is generically exact then we may form
its determinant $D$ (over the function field) in the previous sense, and
the divisor of zeros and poles of the rational function $D$ is the
Div of the complex.
\newline
One reason why this construction is only 'almost' explicit is that there are
choices of minors. Another reason is that it gives a rational
function even in cases, like (4.6) and (5.19), where one knows that the answer
is actually a polynomial.
\newline
\newline
(4.6) Returning to (4.3), suppose that $E$ is a vector bundle of
rank $n+1$ on $X^n$, generated by global sections.
For each $s \in H^0(X,E)$ denote by $K(E,s)$ the Koszul complex,
obtained from multiplication by $s$ in the exterior algebra of $E$
$$
 0 @>>> \Cal O_X @>{\wedge s}>>  E  @>{\wedge s}>>
\wedge^{2} E @>{\wedge s}>> \dots @>{\wedge s}>>
\wedge^{n+1} E @>>> 0
$$
Since a Koszul complex is exact iff it is exact at the 0-th spot,
the resultant $R_E$ is equal to the set of sections $s$
such that  $K(E,s)$ is not exact. For variable $s$, the complex
above may be viewed as a complex $K(E)$ on
$X \times \Bbb P H^0(X,E)$, and denoting by $p$ the projection onto
$\Bbb P H^0(X,E)$, we obtain  a determinantal
expression for the resultant
$$R_E = \text{Div } Rp_* K(E)$$
Since for a complex $F$ and a line bundle $L$
we have $\text{Div } Rp_* F = \text{Div } Rp_* (F \otimes L)$
(\cite{KM}, Prop. 9 (b)), we may
choose a very ample line bundle $\Cal O_X(1)$ on $X$ and
replace $K(E)$ by $K(E) \otimes O_X(m)$
with $m$ large so that $H^i(X, (\wedge^j E)(m)) = 0$
$(i>0, j \ge 0)$. Then $R_E$ is equal to the Div of the complex
(see Cayley-Koszul complex, \cite{GKZ})
$$ 0 @>>> H^0(\Cal O_X(m)) \otimes \Cal O
@>{\wedge s}>>
 H^0(E(m)) \otimes \Cal O(1)
@>{\wedge s}>> \dots @>{\wedge s}>>
 H^0((\wedge^{n+1} E)(m)) \otimes \Cal O(n+1)
@>>> 0$$
on the  projective space $\Bbb P H^0(X,E)$.

\newpage

\flushpar
{\bf \S 5 Hessians.}
\newline
\newline
In the first part of this section, which is independent of the previous
three sections, we prove a result (see (5.11)) on existence and
uniqueness of Hessians.
After this, we set up the determinantal description of Hessians,
which is the main result of this article.
\newline
\newline
(5.1) Let $f: X \to S$ be a smooth and proper morphism of
relative dimension
one, where the varieties $X$ and $S$ are
assumed to be smooth and irreducible.
Let $L$ be a line bundle on $X$ such that $f_*(L)$ is
locally free of rank $n+1$.
Consider the natural homomorphism
$$\tau: f^* f_*(L) \to \Cal P^{n}_{X/S}(L)$$
of locally free modules of rank $n+1$ on $X$ and denote by
$w_{f,L} = \text{det } \tau$ the relative Wronskian
and by $F_{f,L}$ the divisor of zeros of $w_{f,L}$.
\newline
\newline
(5.2) For the linear
equivalence class of $F_{f,L}$ we have
$$
\align
[F_{f,L}] &= c_1\Cal P^{n}_{X/S}(L) - c_1f^* f_*(L) \\
&= (n+1) \ c_1 L + \binom {n+1}2 \ c_1\Omega^1_{X/S} - c_1f^* f_*(L)
\endalign
$$
where the second equality is proved using (2.8).
\newline
\newline
(5.3) Suppose that our family of curves is embeded
in projective space, so that we have a commutative diagram
$$
\CD
X @>i>> S \times \Bbb P^r = Y
 @>q>> \Bbb P^r \\
 @VVfV          @VVpV                   @. \\
S @= S   @.
\endCD
$$
where $i$ is a closed embeding and $p$ and $q$
are the canonical projections. We shall use
the notation
$\Cal O_Y(d) = q^* \Cal O_{\Bbb P^r}(d)$ and
$\Cal O_X(d) = (q \circ i)^* \Cal O_{\Bbb P^r}(d)$
for $d \in \Bbb Z$. The line bundle
$L$ will be chosen to be
$L= \Cal O_X(m)$ for a certain $m \in \Bbb N$.
\newline
\newline
(5.4) If $\Cal I$ is the ideal sheaf of $X$ in $Y$,
from the exact sequence
$$0 \to \Cal I/\Cal I^2 \to
q^*\Omega^1_{\Bbb P^r}|_X = \Omega^1_{Y/S}|_X \to
\Omega^1_{X/S} \to 0 $$
we obtain
$$c_1\Omega^1_{X/S} = (-r-1) \ c_1\Cal O_X(1) -
c_1 \Cal I/\Cal I^2$$
\newline
(5.5) Let us also assume that our family is the family of smooth
complete intersection curves of
multi-degree $(d_1, \dots, d_{r-1})$
in $\Bbb P^r$, so that $S$ is an open subvariety of the product
$P = \prod_{1 \le i \le r-1} P_i$ of projective spaces
$P_i = \Bbb P H^0(\Bbb P^r, \Cal O(d_i))$, and $X$ is given as
$$ X = \{(F_1, \dots, F_{r-1}, x) / F_i(x) = 0 \ \forall i\}
\subset Y$$
We denote
by $\Cal O_i(1)$ the pull-back of $\Cal O_{P_i}(1)$ by the
canonical projection $S \to P_i$. Also, we write
$L_i = p^*\Cal O_i(1) \otimes \Cal O_Y(d_i)$, so that $X$ is
the scheme of zeros of a regular section of the vector bundle
$E = \bigoplus_{i=1}^{r-1} L_i$ on $Y$. It follows that
$\Cal I/\Cal I^2 = E^*|_X$ and hence
$c_1 \Cal I/\Cal I^2 = - c_1 E =
- \sum_i p^*c_1\Cal O_i(1) + d_i c_1\Cal O_Y(1)$. We obtain
from (5.4)
$$c_1\Omega^1_{X/S} =
(-r-1 + \sum_i d_i) \ c_1\Cal O_Y(1) + \sum_i p^*c_1\Cal O_i(1) $$
\newline
(5.6) To determine $c_1 f_*(L)$, we consider the
Koszul resolution
of $\Cal O_X$ tensored by $\Cal O_Y(m)$
$$0 \to  \Cal O_Y(m) \otimes \wedge^{r-1} E^*
\to \dots \to \Cal O_Y(m) \otimes E^*  \to
\Cal O_Y(m) \to L \to 0$$
We start by calculating
$$\wedge^{\alpha} E^* = \wedge^{\alpha} \bigoplus_i L^*_i =
 \bigoplus_{|J|=\alpha}  \bigotimes_{j \in J} L^*_j =
\bigoplus_{|J|=\alpha}
 \Cal O_Y(-d_J) \otimes  p^* \bigotimes_{j \in J} \Cal O_j(-1)$$
where $J$ denotes a subset of $\{1, \dots, r-1\}$ and
$d_J := \sum_{j \in J} d_j$.
It follows from the projection formula that
$$p_*(\Cal O_Y(m) \otimes \wedge^{\alpha} E^*) =
\bigoplus_{|J|=\alpha}
 H^0(\Bbb P^r,\Cal O(m-d_J)) \otimes
 \bigotimes_{j \in J} \Cal O_j(-1)$$
and that
$R^ip_*(\Cal O_Y(m) \otimes \wedge^{\alpha} E^*) = 0$
for $1 \le i \le r-1$. This last vanishing implies that $p_*$
of the resolution of $L$ above is exact, as is easily seen by
breaking the resolution into 3-term exact sequences.
Therefore,
$$
\align
c_1 f_*(L) &=
\sum_{\alpha = 0}^{r-1} (-1)^{\alpha} \
c_1 p_*(\Cal O_Y(m) \otimes \wedge^{\alpha} E^*) \\
&=
 \sum_{\alpha = 0}^{r-1} (-1)^{\alpha}
\sum_{|J| = \alpha} h^0(\Bbb P^r,\Cal O(m-d_J))
\sum_{j \in J} c_1 \Cal O_j(-1) \\
&=
\sum_{\{(J,j) / j \in J \}}
(-1)^{|J|}\  h^0(\Bbb P^r,\Cal O(m-d_J)) \ c_1 \Cal O_j(-1) \\
&=
\sum_{j=1}^{r-1} \left( \shave{\sum_{J \owns j}}
(-1)^{|J|}\  {\eightpoint \binom{m-d_J +r}r} \right)
\ c_1 \Cal O_j(-1)
\endalign
$$
where we use the convention $\binom ab = 0$ if $a < b$.
\newline
\newline
(5.7) It also follows from (5.6) that
$$n+1 = \text{rank }f_*(L) =
\sum_{J \subset \{1, \dots, r-1\}}
(-1)^{|J|} \  {\eightpoint \binom{m-d_J +r}r}$$
\newline
(5.8) Combining the results above, we obtain
$$[F_{f,L}] = a \ c_1 \Cal O_X(1) +
\sum_{j=1}^{r-1} \  b_j \ c_1 \Cal O_j(1)$$
with $a = (n+1)m + \binom {n+1}2 (-r-1+\sum_i d_i)$ as predicted by
the Pl\"ucker formula \cite{ACGH} and
$$b_j = {\eightpoint \binom {n+1}2} + \sum_{J \owns j}
(-1)^{|J|}\  {\eightpoint \binom{m-d_J +r}r} $$
\newline
(5.9) PROPOSITION. Let $\bar F_{f,L}$ and $\bar X$ denote the closures of
$F_{f,L}$ and $X$ in
$P \times \Bbb P^r$, with $P$ as in (5.5). Then the
equality
$$[\bar F_{f,L}] = a \ c_1 \Cal O_{\bar X}(1) +
\sum_{j=1}^{r-1} \  b_j \ c_1 \Cal O_j(1)$$
holds in $\text{Pic}(\bar X)$.
\newline
\newline
PROOF:
Let us recall that $S$ in (5.5)
is the open subvariety of all $F = (F_1, \dots, F_{r-1}) \in P$
such that $F_1= \dots= F_{r-1} = 0$ is a smooth curve, and the point
of the present proof is that now we
have to deal with bad $F$'s that violate this condition.
More precisely, let $\bar f: \bar X \to P$ be the restriction of
the projection $p: \Bbb P^r \times P \to P$ as in (5.3), and denote
by $\Sigma \subset \bar X$ the subvariety where $\bar f$ is not
smooth ($\Sigma$ has codimension at least two in $\bar X$).
 Let $w_{\bar f, L}$ denote the corresponding Wronskian,
defined on $\bar X - \Sigma$, and $F_{\bar f,L}$
its scheme of zeros. The class of
$F_{\bar f,L}$ in
$\text{Pic}(\bar X - \Sigma) \cong \text{Pic}(\bar X)$
is given by the same formula (5.8), since the previous
argument goes through (see \cite{C}).
What we need to check
is that $F_{\bar f,L} - F_{f,L}$ does not have divisorial components,
so that the closures in $\bar X$ of $F_{\bar f,L}$
and $F_{f,L}$ coincide. For this, we only need
to consider phenomena that occur in codimension one in the
smooth variety $\bar X$.
First let us look at the subvariety $J \subset P$
where the dimension of the fibers of $\bar f: \bar X \to P$ jumps: this is
under control since
$\bar f^{-1}(J) \subset \Sigma$ has no divisorial components.
Next, let $\Delta \subset P$ denote the divisor
parametrizing $F$'s such that $F=0$ is a
singular curve, with one node as its only singularities.
 The divisor
$\Delta$ is irreducible, and for
$F \in \Delta$ the corresponding curve $F=0$ is irreducible
and it has points that are not $L$-flexes (i.e. its Wronskian,
defined away from the singular point, is not identically zero).
Therefore, the irreducible divisor $f^{-1}(\Delta)$
is not contained in $\bar F_{f,L}$. Finally, we observe that
$\bar f^{-1}(\bar \Delta - \Delta - J)$ has codimension two, and
 hence the proposition is proved.
\newline
\newline
(5.10) Now we analyse the existence and uniqueness of $m$-Hessians.
We use the notation $\bar Y =  \Bbb P^r \times P$,
$\bar F = \bar F_{f,L}$ and
$\Cal F = \Cal O_{\Bbb P^r}(a) \otimes
\bigotimes_{j=1}^{r-1} \Cal O_j( b_j )$ with $a, b_j$ as in (5.8).
Consider the natural exact sequence of ideal sheaves on $\bar Y$,
twisted by $\Cal F$:
$$ 0 @>>> \Cal F \otimes \Cal I_{\bar X / \bar Y}  @>>>
\Cal F \otimes \Cal I_{\bar F / \bar Y} @>\rho>>
\Cal F \otimes \Cal I_{\bar F / \bar X}
\cong \Cal O_{\bar X} @>>> 0$$
The general linear group $G = GL_{r+1}(k)$ acts on $\bar Y$, preserving
$\bar X$ and $\bar F$, so that the sequence above is of $G$-Modules.
It follows from (5.6) and Kunneth that
$H^1(\bar X, \Cal F \otimes \Cal I_{\bar X / \bar Y})=0$, and taking
cohomology we obtain an exact sequence of finite dimensional
$G$-vector spaces
$$
0 @>>>
\left \{
\aligned & \text{forms vanishing} \\
&\text {on }  \bar X
\endaligned
\right \}
@>>>
\left \{
\aligned & \text{forms vanishing on}  \\
& \text {m-flexes of curves}
\endaligned
\right \}
@>H^0(\rho)>>
H^0(\bar X, \Cal F \otimes \Cal I_{\bar F / \bar X}) =  k. w_{f,L}
@>>> 0
$$
where by forms we mean elements in
$$
\align
H^0(\Bbb P^r \times P, \Cal F) & \cong
H^0(\Bbb P^r, \Cal O(a)) \otimes
\bigotimes_{j=1}^{r-1} \text{Sym}^{b_j}H^0(\Bbb P^r, \Cal O(d_j))^* \\
& \cong
\text{Sym}^a(k^{r+1}) \otimes
\bigotimes_{j=1}^{r-1} \text{Sym}^{b_j} \text{Sym}^{d_j}(k^{r+1})^*
\endalign
$$
An $m$-Hessian is a $G$ semi-invariant element in
$H^0(\bar X, \Cal F \otimes \Cal I_{\bar F / \bar Y})$, and hence
corresponds to a $G$-splitting of the sequence above. Since $G$ is
reductive, such a splitting exists. It also follows that if
$H_m$ denotes an $m$-Hessian
then any other $m$-Hessian may be represented, up to multiplicative
constant, as
$$H_m + \sum_{j=1}^{r-1} G_j \ F_j$$
where $F_j$ is the general
polynomial of degree $d_j$ and $G_j = G_j(x, F_1, \dots, F_{r-1})$
is zero or a $G$ semi-invariant of the appropriate multi-degree, when
they exist.
\newline
\newline
To summarize, we may state:
\newline
\newline
(5.11) PROPOSITION. For any $r, m, d_1, \dots, d_{r-1}$ an $m$-Hessian
$H_m$ exists and is esentially unique. It may be represented as
a form $H_m= H_m(x, F_1, \dots, F_{r-1})$
of multi-degree $(a, b_1, \dots, b_{r-1})$, where $a, b_j$ are as in
(5.8).
\newline
\newline
(5.12) It is verified that the classical formulas in the Introduction
 have the degree stated
in (5.11).
Notice that in all those cases one has $m < d_j$ for all $j$,
and hence $b_j = \binom {n+1}2$.
\newline
\newline
(5.13) Now we turn to the description of the complex of graded free
modules representing an $m$-Hessian.
Let $X \subset \Bbb P^r_k = P$ be a smooth
complete intersection curve defined by the sequence of polynomials
$F_1, \dots, F_{r-1} \in A = k[x_0, \dots, x_r]$
of degrees $d_1, \dots, d_{r-1}$.
Denote $\Cal E = \bigoplus_{i=1}^{r-1} \Cal O_P(d_i)$, so that
$F = (F_1, \dots, F_{r-1}) \in H^0(P, \Cal E)$ and $X = (F=0)$.
The Koszul complex $K_F$
$$0 \to  \Cal O_P(m) \otimes \wedge^{r-1} \Cal E^*
\to \dots \to \Cal O_P(m) \otimes \Cal E^*  \to
\Cal O_P(m)$$
provides a resolution of $\Cal O_X(m)$ by graded free $\Cal O_P$-Modules.
\newline
\newline
(5.14) Let us consider the commutative diagram of complexes of
$\Cal O_X$-Modules
$$
\CD
\Cal H := H^0(P, K_F) \otimes \Cal O_X   @>\rho>> H^0(X,\Cal O_X(m)) \otimes
\Cal O_X @>>> 0 \\
@V{\bar\tau}VV                                        @VV{\tau}V     @.\\
\Cal P := \Cal P^n(K_F)                       @>\rho>>     \Cal P^n_{X/k}(\Cal
O_X(m)) @>>> 0
\endCD
$$
with $\Cal P^n(K_F)$ as in (3.9), using $K_F$ for $(M,\mu)$,
$\bar\tau$ as in (3.10) and $n$ as in (5.7).
\newline
\newline
(5.15) The Euler sequences (2.25) combine to give
 an exact sequence of complexes of $\Cal O_X$-Modules
$$ 0 @>>> \Cal P^n(K_F) @>j>>
\text{P}^n := \iota^* (\text{P}^n_{A/k}(F)\otimes A(m))\sptilde
@>\epsilon>>
\text{P}^{n-1} := \iota^* (\text{P}^{n-1}_{A/k}(F) \otimes A(m)) \sptilde
@>>> 0 $$
where $(\text{P}^n_{A/k}(F)\otimes A(m))\sptilde$ is the complex of (3.5) or
(3.8)
with each module in the complex twisted by $m$ and then sheafified,
 $\epsilon$ is defined by
duals of Euler operators  $.(E-s)^*$ as in (2.27),
for various constants $s$,  and
$\iota: X \to P$ is the inclusion map. More precisely, the exact sequence
above in degree $j$ is
$$
0 @>>> \bigoplus_{|J| = j} \Cal P^{n-j}_{P/k} (\Cal O(m-d_J)) @>\jmath>>
\bigoplus_{|J| = j} \iota^*(\text{P}^{n-j}_{A/k}(A)(m-d_J)\sptilde)
 @>\epsilon_j>>
\bigoplus_{|J| = j} \iota^*(\text{P}^{n-j-1}_{A/k}(A)(m-d_J)\sptilde) @>>> 0
$$
with $\epsilon_j = \bigoplus_{|J| = j}.(E-(m - d_J)^*$, where
$d_J = \sum_{j \in J} d_j$ and $J$ is a subset of $\{1, \dots, r-1\}$.
\newline
\newline
(5.16) REMARK. The complexes $(\text{P}^n_{A/k}(F)\otimes A(m))\sptilde$
are defined on $P$, and the maps $\epsilon$ are defined on $P$,
but they commute with the differentials of the complexes only after
applying $\iota^*$
(i.e. modulo multiples of the $F_i$'s).
\newline
\newline
(5.17) It is easy to check that the upper row of (5.14) is exact.
The lower row is also exact, as it follows from (3.6)
and the long exact sequence of homology
for (5.15). In other words,
the maps $\rho$ are quasi-isomorphisms, and hence
$\text{Div}(\rho) = 0$. Now \cite{KM} (Thm. 3 (i), p. 48) implies that
$$ \text{Div }(\tau) = \text{Div }(\bar\tau)$$
which achieves the first step in the direction of lifting
the Wronskian $\text{det}(\tau)$ to $\Bbb P^r$.
\newline
\newline
(5.18) Next we would like to replace
$\bar\tau: \Cal H \to \Cal P$ by a homomorphism of
complexes of graded free modules.
This will be done by a cone construction. Consider the diagram
$$
\CD
       @.      \Cal H        @.                     @.
@.   \\
 @.        @V{\bar\tau}VV            @.                      @.
   @. \\
0     @>>>   \Cal P      @>j>>       \text{P}^n    @>\epsilon>>
\text{P}^{n-1}     @>>>    0  \\
@.        @V{\alpha}VV            @.                      @.
@. \\
@.  \Cal C
\endCD
$$
where $\Cal C = \text{cone}(\epsilon)(1) = \text{P}^n \oplus \text{P}^{n-1}(1)$
is the cone complex of $\epsilon$
shifted by one. It follows from \cite{Bo,2} (Ex. 9, p. 174) that
the natural map
$\alpha(x) = (j(x), 0)$ is a quasi-isomorphism,
and therefore $\text{Div}(\alpha) = 0$.
Defining $\tilde{\tau} = \alpha \circ \bar\tau$  we have the equality
$$\text{Div}(\bar\tau) = \text{Div}(\tilde\tau)$$
We state the result in the following way
\newline
\newline
(5.19) THEOREM. If $\tilde{\tau}$
denotes the natural homomorphism of complexes
of $\Cal O_X$-Modules
$$H^0(P, K_F) \otimes \Cal O_X @>{\tilde{\tau}}>>
\iota^*(\text{P}^n_{A/k}(F) \otimes A(m))\sptilde   \oplus
\iota^*((\text{P}^{n-1}_{A/k}(F) \otimes A(m)) \sptilde )(1)$$
then $\text{Div}(\tilde{\tau})$ is a rational $m$-Hessian for $X$.
\newline
\newline
We need to explain what we mean by rational $m$-Hessian and why
they appear in our statement.
As mentioned in (5.16), the differentials
of the complexes in question, and the map $\tilde{\tau}$, are all
 given by matrices of homogenous polynomials, that is, they are
defined over $P$. The difficulty is that our cone $\Cal C$ is a complex
only on $X$, that is, the condition $d^2=0$ is satisfied on $P$ only modulo
the $F_i$'s. To calculate $\text{Div}(\tilde{\tau})$ we  carry
out the process of (4.5) with these matrices of polynomials. To do this
we need to make a choice of admissible minors. We don't need the
condition $d^2=0$ to calculate the alternating product of determinants
of the chosen minors, but different choices will give different answers,
all equal on $X$ (by (5.18), compare (5.11)).
These various "Div's" are
quotients $A/B$ of homogeneous polynomials on $P$.
We obtain a finite number of rational $m$-Hessians and
our computations seem to indicate that none of these is a
polynomial Hessian. It would be desirable to find a further manipulation
of these complexes that leads to a polynomial Hessian like in (5.11).
\newline
\newline
(5.20) For plane curves the diagrams above take a simpler form
that perhaps it is worth making explicit. Let
$X = (F=0) \subset P = \Bbb P^2$ with $F$ of degree $d$ and fix
$m \in  \Bbb N$. The basic diagram is
$$
\minCDarrowwidth{.4 cm}
\eightpoint
\CD
0 @>>> H^0(P, \Cal O(m-d)) \otimes \Cal O_X
@>{.F}>> H^0(P, \Cal O(m)) \otimes \Cal O_X
@>>> H^0(X, \Cal O(m)) \otimes \Cal O_X
@>>> 0      \\
@.  @V{\tau}VV         @V{\tau}VV  @V{\tau}VV    @.    \\
0 @>>> \iota^*\Cal P^{n-1}_{P/k}\Cal O(m-d)  @>>>
\iota^*\Cal P^{n}_{P/k}\Cal O(m)  @>>> \Cal P^{n}_{X/k}\Cal O(m) @>>> 0 \\
@. @VVV         @VVV          \\
 @. \iota^*\text{P}^{n-1}_{A/k}(m-d)\sptilde  @>{.\Delta F}>>
\iota^*\text{P}^{n}_{A/k}(m)\sptilde  \\
 @.  @V{\epsilon}VV         @V{\epsilon}VV    \\
 @. \iota^*\text{P}^{n-2}_{A/k}(m-d)\sptilde  @>{.\Delta F}>>
\iota^*\text{P}^{n-1}_{A/k}(m)\sptilde  \\
@. @VVV         @VVV        \\
 @. 0  @. 0
\endCD
$$
where $.\Delta F$ is multiplication by
$\Delta F = \sum_{\alpha} \partial^{\alpha}(F) \  (dx)^{\alpha}$
in the symmetric algebra of $\Omega^1_{A/k}$
(see (2.34) or (3.8)). According to  our previous analysis, the divisor
of $m$-flexes is the Div of the double complex
$$
\eightpoint
\CD
H^0(P, \Cal O(m-d)) \otimes \Cal O_X
@>{.F}>> H^0(P, \Cal O(m)) \otimes \Cal O_X \\
 @V{\tau}VV         @V{\tau}VV   \\
 \iota^*\text{P}^{n-1}_{A/k}(m-d)\sptilde  @>{.\Delta F}>>
\iota^*\text{P}^{n}_{A/k}(m)\sptilde  \\
 @V{\epsilon}VV         @V{\epsilon}VV   \\
 \iota^*\text{P}^{n-2}_{A/k}(m-d)\sptilde  @>{.\Delta F}>>
\iota^*\text{P}^{n-1}_{A/k}(m)\sptilde \\
\endCD
$$
that is, the Div of the associated total complex
$$
\align
H^0(P, \Cal O(m-d)) \otimes \Cal O_X @>{\delta_2}>> &
H^0(P, \Cal O(m)) \otimes \Cal O_X \ \oplus \
\iota^*\text{P}^{n-1}_{A/k}(m-d)\sptilde
@>{\delta_1}>> \\
& \iota^*\text{P}^{n}_{A/k}(m)\sptilde \  \oplus \
\iota^*\text{P}^{n-2}_{A/k}(m-d)\sptilde
@>{\delta_0}>> \iota^*\text{P}^{n-1}_{A/k}(m)\sptilde
\endalign
$$
where $\delta_2(a) = (F.a, - \tau(a))$,
$\delta_1(b,c) = (c.\Delta F + \tau(b), -\epsilon(c))$ and
$\delta_0(d,e) = e.\Delta F + \epsilon(d)$.


\newpage

\Refs
\widestnumber \key {A-C-G-H,2}
\midspace{0.1 in}

\ref
\key {\bf ACGH}
\by E. Arbarello, M. Cornalba, P. Griffiths and J. Harris
\paper Geometry of Algebraic Curves
\jour Springer-Verlag
\yr 1984
\endref

\ref
\key {\bf BGS}
\by J.-M. Bismut, H. Gillet and C. Soule
\paper Analytic torsion and holomorphic determinant bundles
\jour Comm. Math. Physics
\yr 1988
\endref

\ref
\key {\bf BB}
\by W. Borho and J.-L. Brylinski
\paper Differential operators on homogeneous spaces
\jour Inventiones Mathematicae
\vol 69
\yr 1982
\endref

\ref
\key {\bf B}
\by A. Borel et al.
\paper Algebraic D-modules
\jour Academic Press
\yr 1987
\endref

\ref
\key {\bf Bo,1}
\by N. Bourbaki
\paper Algebra, Chapter II
\jour Springer-Verlag
\yr 1989
\endref

\ref
\key {\bf Bo,2}
\by N. Bourbaki
\paper Algebre, Chapitre X
\jour Masson
\yr 1980
\endref

\ref
\key{\bf Bu}
\by I. Bucur
\paper Triangulated categories and algebraic K-theory
\jour  Springer Lecture Notes
\vol 108
\endref

\ref
\key{\bf Ca}
\by A. Cayley
\paper On the sextatic points of a plane curve
\vol V (341)
\jour Collected Mathematical Papers
\endref

\ref
\key{\bf C}
\by F. Cukierman
\paper Families of Weierstrass points
\jour Duke Math. J.
\yr 1989
\endref

\ref
\key{\bf Dy,1}
\by R. H. Dye
\paper Osculating primes to curves of intersection
in 4-space, and to certain curves in n-space
\vol 18
\yr 1973
\jour Proc. Edinburgh Math. Soc.
\endref

\ref
\key{\bf Dy,2}
\by R. H. Dye
\paper The hyperosculating spaces to certain curves in [n]
\vol 19
\yr 1974/75
\jour Proc. Edinburgh Math. Soc.
\endref

\ref
\key{\bf E,1}
\by W. L. Edge
\paper The osculating solid of a certain curve in [4]
\vol 17
\yr 1971
\jour Proc. Edinburgh Math. Soc.
\endref

\ref
\key{\bf E,2}
\by W. L. Edge
\paper The osculating spaces of a certain curve in [n]
\vol 19
\yr  1974/75
\jour Proc. Edinburgh Math. Soc.
\endref

\ref
\key{\bf GKZ,1}
\by I. Gelfand, M. Kapranov and A. Zelevinsky
\paper Projectively dual varieties and hyperdeterminants
\yr 1989
\vol 39 (2)
\jour Soviet Math., Doklady
\endref

\ref
\key{\bf GKZ,2}
\by I. Gelfand, M. Kapranov and A. Zelevinsky
\paper Hyperdeterminants
\yr 1992
\jour Advances in Math.
\endref

\ref
\key{\bf G}
\by A. Grothendieck
\paper Elements de calcul infinitesimal
\jour Seminaire Cartan
\yr 1960/61
\endref

\ref
\key {\bf EGA IV}
\by A. Grothendieck and J. Dieudonne
\paper EGA IV (4)
\jour Publ. Math. IHES
\vol 32
\yr 1967
\endref

\ref
\key {\bf I}
\by L. Illusie
\paper Complexe Cotangent et Deformations I
\jour  Springer Lecture Notes
\vol 239
\endref

\ref
\key{\bf KM}
\by F. Knudsen and D. Mumford
\paper The projectivity of the moduli space of stable
 curves: Preliminaries on det and div
\yr  1976
\jour Math. Scandinavica
\endref

\ref
\key{\bf Mi}
\by J. Milnor
\paper Whitehead torsion
\yr 1966
\jour Bulletin A.M.S
\endref

\ref
\key{\bf P}
\by R. Piene
\paper Numerical characters of a curve in projective n-space
\yr 1976
\jour Proc. Nordic Summer School, Oslo
\endref

\ref
\key{\bf R}
\by A. Ranicki
\paper The algebraic theory of torsion I. Foundations
\vol 1126
\jour  Springer Lecture Notes
\endref

\ref
\key{\bf S}
\by G. Salmon
\paper Analytic geometry of three dimensions
\jour Chelsea
\endref

\ref
\key{\bf W}
\by R. Walker
\paper Algebraic Curves
\jour Dover Publications
\endref

\endRefs
\enddocument